\documentclass[preprint,showpacs,showkeys,preprintnumbers,amsmath,amssymb]{revtex4}

\usepackage{graphicx}
\usepackage{dcolumn}
\usepackage{bm}

\begin{document}

\preprint{}

\title{A Numerical Examination of the Castro-Mahecha Supersymmetric 
Model of the Riemann Zeros}

\author{Paul B. Slater}%
\email{slater@kitp.ucsb.edu}
\affiliation{%
ISBER, University of California, Santa Barbara, CA 93106\\
}%
\date{\today}

\begin{abstract}
The unknown parameters of the recently-proposed 
({\it Int J. 
Geom. Meth. Mod. Phys.} {\bf{1}}, 751 [2004])
Castro-Mahecha model of 
the imaginary parts ($\lambda_{j}$) 
of the nontrivial Riemann zeros are the phases 
($\alpha_{k}$) and the frequency parameter ($\gamma$) 
of the Weierstrass function of {\it fractal} 
dimension $D=\frac{3}{2}$ and the turning points ($x_{j}$) 
of the {\it supersymmetric} potential-squared $\Phi^2(x)$ --- which incorporates the
{\it smooth} 
Wu-Sprung potential ({\it Phys. Rev.} E {\bf{48}}, 2595 [1993]), giving the {\it average} level density of 
the Riemann zeros. We conduct numerical investigations
to estimate/determine these parameters --- as well as a parameter 
($\sigma$) we introduce to {\it scale} the fractal contribution.
Our primary analyses involve 
two sets of {\it coupled} equations: one 
set being of the form 
$\Phi^{2}(x_{j}) = \lambda_{j}$, and the other set corresponding to 
the 
fractal extension --- according to an {\it ansatz} of 
Castro and Mahecha -- of 
the Comtet-Bandrauk-Campbell 
(CBC) quasi-classical quantization conditions for good supersymmetry, 
$\frac{2}{\Gamma(D/2)}
\int_{-x_{j}}^{x_{j}} 
\frac{[\lambda_{j} - \Phi^{2}(x')]^{1/2}}{(x_{j}-x')^{1-D/2}} d x' = j  \pi$. 
Our analyses suggest the possibility strongly that $\gamma$ converges to
its theoretical lower bound of 1, and the possibility that all the phases
($\alpha_{k}$) 
should be set to zero. 
We also uncover interesting formulas for certain fractal turning points.
\end{abstract}


\maketitle

\section{Introduction}
The precise nature of the {\it zeros} of the Riemann zeta function is of paramount
mathematical \cite{conrey}, as well as physical interest \cite{circular}.
Of course, perhaps the most famous, not yet fully 
resolved,  mathematical conjecture
(the Riemann Hypothesis \cite{rockmore}) is that {\it all} 
the nontrivial zeros ($\eta_{j} + i \lambda_{j}$) --- the trivial zeros 
simply equalling $-2 j$, $j=1,2,3,\ldots$ --- have 
real components $\eta_{j}= \frac{1}{2}$.

Castro and Mahecha \cite{carlos1} recently proposed a 
{\it supersymmetric}  implementation of the 
Wu-Sprung model \cite{wu} of the $\lambda_{j}$'s, employing
a Weierstrass (continuous and nowhere differentiable) 
function \cite{berry} for the {\it fractal} structure of
dimension $D=\frac{3}{2}$  that Wu and Sprung had found for their 
local one-dimensional potential ($V$) (cf. \cite{khuri,rosu}).
The Wu-Sprung potential $V$ --- which generates the smooth {\it average} 
level density obeyed by the Riemann zeros --- satisfied 
{\it Abel's integral equation} 
\cite[eq. (6)]{wu}, and was
written implicitly as \cite[eq. (7)]{wu},
\begin{equation} \label{WSpotential}
x=\frac{1}{\pi} \Big( \sqrt{V-V_{0}} \ln \frac{V_{0}}{2 \pi e^2} +  \sqrt{V}
\ln \frac{\sqrt{V}+\sqrt{V-V_{0}}}{\sqrt{V}-\sqrt{V-V_{0}}} \Big).
\end{equation}
Here $V_{0}= 3.10073 \pi \approx 9.74123$.

The Weierstrass fractal function employed  by Castro and Mahecha
\cite[eq. (81)]{carlos1}, 
\begin{equation} \label{Weierfractal}
W(x,\gamma,D,\alpha_{k})= \Sigma_{k=0}^{\infty}\frac{1- e^{i x \gamma^k}}{\gamma^{k (2 -D)}} e^{2 \pi i \alpha_{k}},
\end{equation}
has both an unknown (countably) infinite set of phases $2 \pi \alpha_{k}$ and
an unknown parameter $\gamma$, which determines the frequencies $\gamma^{k}$.
Further, $D$ is the fractal dimension, which Castro and Mahecha --- following 
the  box-counting argument of 
Wu and Sprung \cite{wu}
(cf. \cite{vanZyl}) --- took to be $\frac{3}{2}$. (Castro and Mahecha also
noted \cite[sec. 6]{carlos1} that this choice of $D=\frac{3}{2}$ 
does,  appealingly,  
yield [the omnipresent] $\frac{1}{f}$ noise.) The variable $x$ parameterizes
the {\it inverted} Wu-Sprung potential $V_{WS}(x)$, 
which in accordance with the 
supersymmetric requirement that the energy be zero in the ground state, 
was translated to be zero at $x=0$. 

Rather than {\it directly} addressing the formidable problem  
``of factoring the ordinary Schr\"odinger equation studied 
by Wu-Sprung'', Castro and Mahecha (CM) postulated the model 
\cite[eq. (82)]{carlos1},
\begin{equation} \label{SUSY}
\Phi^2(x)=V_{WS}(x) +\frac{1}{2} [W(x,D,\gamma,\alpha_{k}) +
W(-x,D,\gamma,\alpha_{k}) + c. c.] +\phi_{0},
\end{equation}
with
\begin{equation}
\Phi^2(x_{j}) = \lambda_{j}.
\end{equation}
($c. c.$ denotes complex conjugation.) 
The additive constant $\phi_{0}=-V_{0}$
is chosen so that, in accordance with 
supersymmetric principles, $\Phi^2(0)=0$. 
The ``turning points'' $x_{j}$ (where the potential energy is
 equal to the total energy), along with the phases $\alpha_{k}$ and the 
parameter $\gamma$, are the unknown (countably infinite) 
parameters of the CM model, 
which we hope here to estimate/determine, by satisfying the equations 
(\ref{SUSY}), {\it simultaneously} with the {\it quantization} conditions 
\cite[eq. (87)]{carlos1},
\begin{equation} \label{CBC}
I_{j}(x_{j},\lambda_{j}) \equiv 2 \frac{1}{\Gamma(\beta)}
\int_{-x_{j}}^{x_{j}} \frac{[\lambda_{j} - \Phi^{2}(x')]^{1/2}}{(x_{j}-x')^{1-\beta}} d x' = j  \pi .
\end{equation}
Castro and Mahecha proposed (cf. \cite[eq. (36)]{laskin}) 
that these relations would constitute 
the {\it fractal} extension of the fermionic phase path integral
approximation (the CBC [Comtet-Bandrauk-Campbell] 
formula (cf. \cite[eq. (1)]{inomata} \cite[eq. (5)]{wu})) --- the 
supersymmetric counterpart of the well-known {\it WKB} quantization rule 
\cite[eq. (2)]{inomata}.
Here $\beta=\frac{D}{2} =\frac{3}{4}$ and $j$ runs over the positive integers.

The Castro-Mahecha model (\ref{SUSY}) was cast in a 
fractal supersymmetric quantum-mechanical (SUSY-QM) 
setting, in an 
effort to 
implement the {\it Hilbert-Polya} proposal \cite{edwards} 
to prove the Riemann Hypothesis, 
under which an Hermitian operator, hypothetically, can 
be found to reproduce the $\lambda_{j}$'s as  its {\it spectrum}.
It had appeared to CM
quite difficult to 
directly solve the SUSY Schr\"odinger equation \cite[eq. (84)]{carlos1},
\begin{equation} \label{difficult}
(\mathcal{D}^{(\beta)}+\Phi) 
(-\mathcal{D}^{(\beta)} + \Phi) \psi_{j}^{(+)}(x) 
= \lambda_{j} \psi_{j}^{(+)}(x),
\end{equation}
where $\hbar= 2 m =1$, $\lambda_{j}^{(+)} = \lambda_{j}^{(-)}$, 
the $\psi_{j}^{(+)}$ are eigenfunctions, 
and the {\it fractional} derivative 
\begin{equation}
\mathcal{D}^{(\beta)}F(t) = \frac{1}{\Gamma(1-\beta)} \frac{d}{dt} 
\int_{-\infty}^t \frac{F(t')}{(t-t')^\beta} dt',
\end{equation}
is intended (cf. \cite{laskin}). (M. Trott suggested that one might {\it directly} pursue the solution of this system, most effectively, 
using a {\it Gr\"unwald-Letnikov} discrete approximation 
\cite{mahori}.)
CM, therefore, 
had recourse to the CBC formula (\ref{CBC}) as an {\it approximation} to the 
implementation of the quantum
{\it inverse scattering} method, which would yield the potential $\Phi$.
Given an exact potential, one could construct the SUSY 
Schr\"odinger equation, giving the $\lambda_{j}$'s for its eigenvalues.
This would realize --- CM claimed --- the 
Hilbert-Polya proposal to prove the Riemann
Hypothesis.

CM reasoned 
 that if the potential energy has the shape of a fractal curve,
then it is very natural to assume that the kinetic energy operator should
be a fractal; namely it should be expressed in terms of {\it fractal} derivatives.

In sec.~\ref{Preliminary}, we begin our numerical analyses by studying
the solutions to (\ref{SUSY}). Then, in sec.~\ref{simultaneous}, we
extend the analyses to include (simultaneously) the two sets of 
equations (\ref{SUSY}) and (\ref{CBC}). 
The {\it fractal} Wu-Sprung potential $V$ is examined in sec.~\ref{potential}, 
in particular its turning points, for a number of which interesting
formulas are obtained. 
In
sec.~\ref{Scale}, we expand the CM model to include a certain scaling
factor. A discussion of our analyses is contained in sec.~\ref{discussion}.
We also include an Appendix outlining recent work of D. Dominici concerning
the {\it inversion} of the Wu-Sprung potential \cite{diego}.
\section{Preliminary analyses} \label{Preliminary}
\subsection{Turning points $x_{j}^{(1)}$ --- for smooth 
potential --- held fixed}
Our initial numerical 
analyses  of the CM  model 
consisted of, for various integral cut-offs $m$ 
on the running index $k$ 
in the Weierstrass fractal function summation (\ref{Weierfractal}),
finding an optimizing set of $\alpha_{k}$'s, while keeping the
$x_{j}$'s {\it fixed} and equal to ($n$ of) the turning points for the 
translated inverted {\it smooth}
Wu-Sprung potential, that is, those $x_{j}^{(1)}$'s for which 
(cf. (\ref{SUSY})),
\begin{equation} \label{start}
V_{WS}(x_{j}^{(1)}) +\phi_{0} =\lambda_{j}. 
\end{equation}
The set of estimated $\alpha_{k}$'s and $\gamma$
would be chosen, so as 
to minimize the sum of squared deviations --- since it did not appear 
possible to achieve complete equality --- between the 
LHS (left-hand side) and RHS of (\ref{start}). Computationally-speaking 
we employed the NMinimize command
of Mathematica, using the ``DifferentialEvolution'' and ``InitialPoints''
options, and for the initial points only {\it increasing} sets of values 
for the $\alpha_{k}$'s, in order to reduce search times. 
Also, the {\it necessary} inversion of the Wu-Sprung potential 
(\ref{WSpotential}) --- in order to 
obtain $V_{WS}(x)$ --- was performed using
the 
Interpolation command of Mathematica (cf. Appendix).

We, first, found certain evidence, at least in this initial context, that
suggested taking 
$\gamma=3 = 2 D$. (In \cite{carlos1}, CM  wrote that 
``it would be intriguing to see if $\gamma=1 +\phi \approx 1.618$,
the inverse of the golden mean, since the golden mean appears in the theory
of quantum noise related to the Riemann Hypothesis'' 
(cf. \cite[p. 481]{berry}).

\begin{figure}
\includegraphics{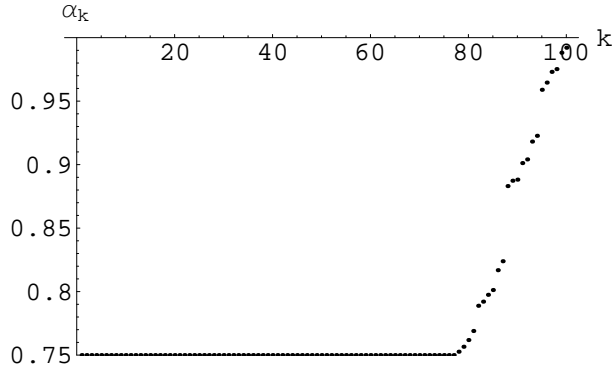}
\caption{\label{fig:100}The one hundred ordered 
phases (divided by $2 \pi$) of the Weierstrass fractal function 
(\ref{Weierfractal})
obtained by minimizing the sum-of-squares deviation from the 
imaginary parts of the first 
one hundred nontrivial Riemann zeros
of the Castro-Mahecha supersymmetric model (\ref{SUSY})}
\end{figure}

\begin{figure}
\includegraphics{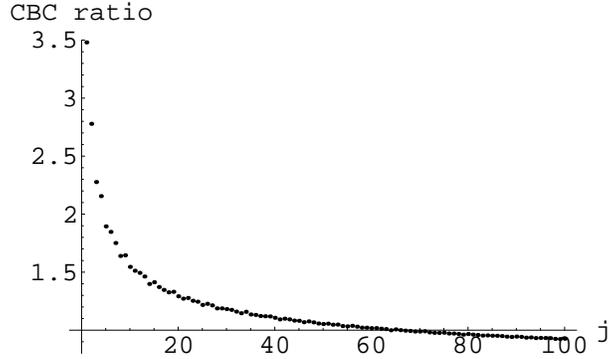}
\caption{\label{fig:quantization}Ratios of the 
(real parts of the) computed values of the
left-hand side of the fractal CBC quantization conditions (\ref{CBC}) --- based on 
the (fixed) turning points, those $x_{j}$ for which $\lambda_{j} 
= V_{WS}(x_{j})$, of the {\it smooth} inverted 
(translated) Wu-Sprung potential and the scaled phases shown in 
Fig.~\ref{fig:100} --- to the 
right-hand side of (\ref{CBC}) (that is, $j \pi$). The ratios would all be 
unity if the
conditions were fully satisfied.}
\end{figure}

Initially, we tried the case $n=100, m=15$, 
insisting that
the phases be ordered from lowest to highest.
Our results were encouraging, and we proceeded to 
successively higher integral values
of $m$,  reaching $m=100$. 
(We make the obvious comment that here the number $m$ of unknowns [$\alpha_{k}$]
{\it equals} the number of [nonlinear] equations.)
For that case, the 
best fit to the first 100 $\lambda_{j}$'s yielded a 
{\it remarkably} small sum-of-squares deviations of 
$2.68927 \cdot 10^{-14}$. 
The first {\it seventy-seven} of the hundred 
ordered phases (that is, $2 \pi \alpha_{k}, k=1,\ldots,70$) 
were {\it all} equal to $\frac{3 \pi}{2} = 2 \pi \beta $ to, typically,  
seven decimal places
(Fig.~\ref{fig:100}).
When we inserted the turning points ($x_{j}^{(1)}$) derived 
from the (smooth) 
inverted 
Wu-Sprung potential for the bounds of integration  --- that 
is, those satisfying (\ref{start}) --- into the fractal CBC
formula (\ref{CBC}), we found that these 
quantization relations  were 
clearly {\it not} fully met. (In general, the 
CBC quantization conditions appear
to hold for a number of ``shape-invariant'' potentials \cite{inomata} 
(cf. \cite{bhaduri}).)
In Fig.~\ref{fig:quantization}, we plot the (CBC) {\it ratio}
of the computed value of the real part of the 
LHS of (\ref{CBC}) to its
RHS (that is, $j \pi$). We see that for the lower values 
of $j$, the LHS's {\it exceed} $j \pi$, while for the higher values, 
they are less than $j \pi$ --- all in an apparently smooth 
 monotonically decreasing manner.

\begin{figure}
\includegraphics{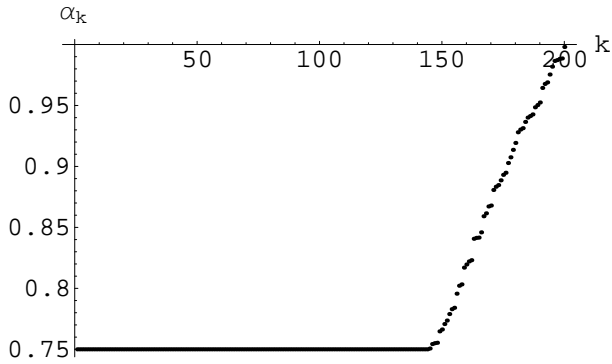}
\caption{\label{fig:200}The {\it two} hundred ordered 
phases (divided by $2 \pi$) of the Weierstrass fractal function
(\ref{Weierfractal})
obtained by minimizing the sum-of-squares deviation from the
imaginary parts of the first
two hundred nontrivial Riemann zeros
of the Castro-Mahecha supersymmetric model (\ref{SUSY})}
\end{figure}

\begin{figure}
\includegraphics{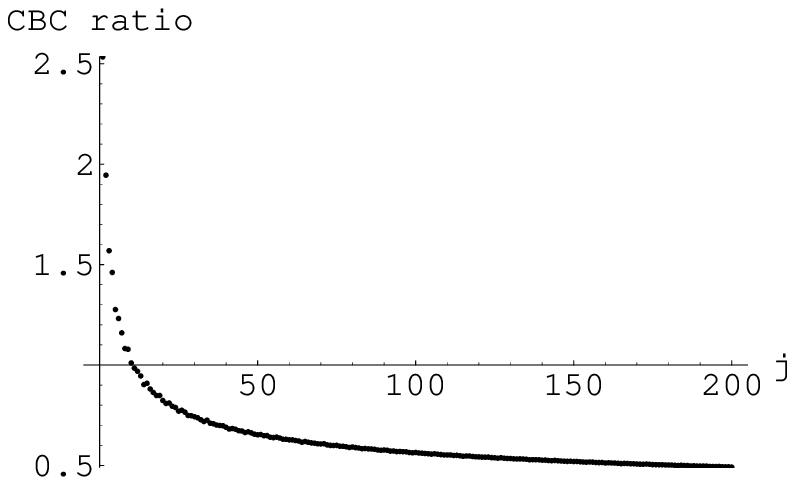}
\caption{\label{fig:quantization2}Ratios, for the case $n=200$, of the
(real parts of the) computed values of the
left-hand side of the fractal 
CBC quantization conditions (\ref{CBC}) to the
right-hand side of (\ref{CBC}) (that is, $j \pi$). The ratios would all be 
unity if the
conditions were fully satisfied.}
\end{figure}

We, then, were able to extend our 
computer analyses to the case $n= m =200$,
with very analogous results
(Fig.~\ref{fig:200}). (The sum of squared deviations was quite small, as 
in the $n=m=100$ analysis, 
$3.48792 \cdot 10^{-14}$.) Most of the 
estimated values of $\alpha_{k}$ were once again
$\frac{3}{4}$ to high precision, this time for $k=1,\ldots,144$.
In Fig.~\ref{fig:quantization2} we present the analogue of Fig.~\ref{fig:quantization} for the $n= m =200$ case.

In regard to the phases of the Weierstrass fractal function,
Castro and Mahecha \cite[sec. 5]{planat} 
had written that it ``is warranted to see
if the statistical distribution of these phases $\alpha_{k}$ 
has any bearing to random matrix theory (the circular unitary
random matrix ensemble) and the recent studies of
quantum phase-locking, entanglement, Ramanujan sums and cyclotomy studied
by \cite{planat}''.

\begin{figure}
\includegraphics{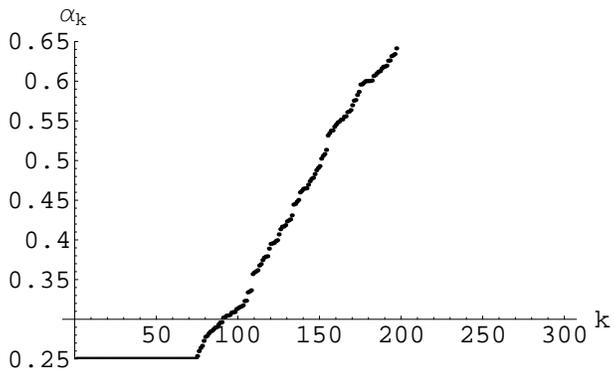}
\caption{\label{fig:300}The {\it three} hundred ordered
phases (divided by $2 \pi$) of the Weierstrass fractal function
(\ref{Weierfractal})
obtained by minimizing the sum-of-squares deviation from the
imaginary parts of the first
three hundred nontrivial Riemann zeros
of the Castro-Mahecha supersymmetric model (\ref{SUSY})}
\end{figure}

When we proceeded to the still higher-order case $n=300,m=300$, we found
superficially, at least, different types of results.
In six different runs (using different random seeds) 
of our Mathematica program, we obtained sum-of-squares deviations
no greater than 
\newline $1.79962 \cdot 10^{-13}$, all six yielding essentially the
same-looking plot. The smallest sum-of-squares was $1.62071 \cdot 10^{-13}$,
corresponding to Fig.~\ref{fig:300}.
The first seventy-four scaled phases ($\alpha_{k}$) were all equal to
$\frac{1}{4}$ (rather than $\frac{3}{4}$, as previously) 
to high accuracy.
Let us note, in any case though, that for 
{\it either} $\alpha_{k}=\frac{1}{4}$ or
$\frac{3}{4}$, the factors $e^{2 \pi i \alpha_{k}}$ in the expansion 
(\ref{Weierfractal})
of the Weierstrass fractal function have {\it zero} real part.

\subsection{Comment}
What it appeared, in retrospect, we had
actually principally 
accomplished above (Figs.~\ref{fig:100}, \ref{fig:200}, \ref{fig:300}) 
was to find sets of phases ($2 \pi \alpha_{k}$) for various $n$ (100, 200 
and 300), 
such 
that the real part of the Weierstrass fractal function (\ref{Weierfractal})
is very close to zero  --- as revealed by plots --- for {\it any} $x$ 
(or $\gamma>1$), be it a turning point or not.
So, these results immediately pertain to the Weierstrass fractal function
and not the Riemann zeros themselves. Nevertheless, 
our estimated phases ($\alpha_{k}$'s) can be employed as 
the {\it first} step of an {\it iterative} procedure --- as discussed 
in the next section.

Since the Castro-Mahecha supersymmetric 
model of the Riemann zeros
requires the {\it simultaneous} obtaining/fitting of phases ($\alpha_{k}$'s), the 
parameter $\gamma$,  
and turning points ($x_{j}$'s), we turn to this problem now.

\section{{\it Joint} Estimation of Phases ($\alpha_{k}$) and Turning Points 
($x_{j}$)} \label{simultaneous}
\subsection{Iterative approach}
\subsubsection{$n=m=100$ analyses}
We would like now to regard the {\it two} sets of
equations (\ref{SUSY}) --- the {\it one} set we have tried to 
satisfy in the above analyses --- {\it and} (\ref{CBC}) 
as a {\it coupled} system and attempt to {\it simultaneously} solve for the 
phases ($\alpha_{k}$) {\it and} turning points ($x_{j}$) 
(still holding $\gamma=3$).
In this context, we again considered the case 
$n=100,m=100$, and took the solution portrayed in Fig.~\ref{fig:100} as the
first step of an {\it iterative} procedure. Using this set of
$\alpha_{k}^{(1)}$'s, we searched for a {\it new} set of turning points
($x_{j}^{(2)}, j=1,\ldots,100$), such that the CBC 
quantization conditions (\ref{CBC}) would {\it all} be exactly 
(up to our numerical accuracy)
satisfied. Then, with this revised set of $x_{j}$'s, we hoped to find  anew the 
set of $\alpha_{k}^{(2)}$'s satisfying (\ref{SUSY}). By reiterating this 
two-step procedure a sufficiently large number of times, we anticipated possibly arriving 
at a {\it final} 
set of $\alpha_{k}^{(N)}$'s and $x_{j}^{(N)}$'s  {\it simultaneously}
satisfying (\ref{SUSY}) and (\ref{CBC}) (cf. \cite{ludger}).
A test for convergence to such a final solution 
would be that the two sets of parameters would only 
change negligibly upon further iterations. \footnote{M. Trott suggested --- in 
analogy to the self-consistent
solution of Schr\"odinger and Poisson equations for some charged 
system --- that it might
be a more superior/stable iterative scheme, to at each stage {\it combine}
the current, 
revised estimates of the $\alpha_{k}^{(i)}$'s and $x_{j}^{(i)}$'s 
with those from
the preceding step of the iterative procedure, in fact, perhaps much more
strongly weighting the preceding ($i-1$) estimates. We have not yet 
thoroughly pursued this possibility, though.}

\begin{figure}
\includegraphics{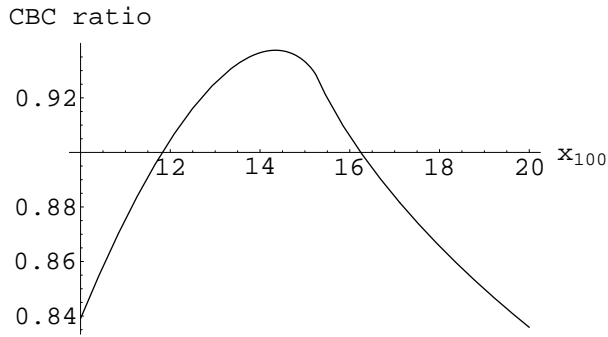}
\caption{\label{fig:Peak100}CBC ratio --- incorporating the phases 
shown in Fig.~\ref{fig:100} --- as a function of the hundreth turning
point ($x_{100}$). The original value of $x_{100}^{(1)}$, 
giving a ratio of 0.926293,  
was 15.315. 
This can be (only) slightly improved to 0.937452, 
by choosing $x_{100}^{(2)}=14.3452$.}
\end{figure}

At the very beginning of the iterative procedure we have 
just outlined, 
we sought to obtain a (hypothetical, second)  set of one hundred 
turning points ($x_{j}^{(2)}$), which would render the CBC quantization
conditions (\ref{CBC}) {\it completely} satisfied when the first set of phases 
($2 \pi 
\alpha_{k}^{(1)}$) --- those depicted (in their scaled form) 
 in Fig.~\ref{fig:100} --- were employed.
However, this proved to be not totally doable.
For example, for the original value of $x^{(1)}_{100}= 15.315$, 
the corresponding
CBC ratio was 0.926293. But when we tried to adjust $x^{(1)}_{100}$ to increase this
ratio closer to the ideal value of 1,  we were only able to obtain a 
{\it slight}
improvement by going to $x^{(2)}_{100} =14.3452$, for which the ratio was 
0.937452.
We show this phenomenon in Fig.~\ref{fig:Peak100}.

\begin{figure}
\includegraphics{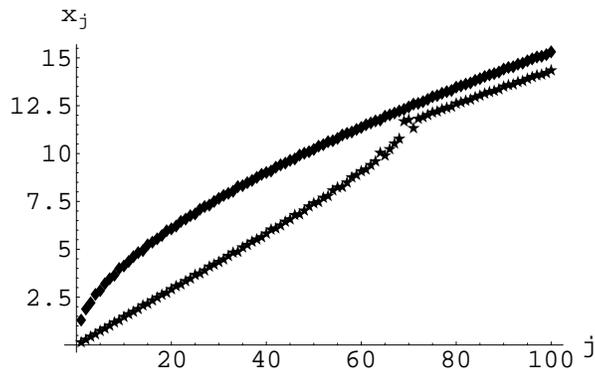}
\caption{\label{fig:Improvement}The original 100 
turning points ($x^{(1)}_{j}$),
forming the higher/dominant line, and the 100 adjusted/revised 
 turning points ($x_{j}^{(2)}$)
that yield CBC ratios optimally closer to {\it unity} --- {\it given} 
the phases 
shown in Fig.~\ref{fig:100}}
\end{figure}

\begin{figure}
\includegraphics{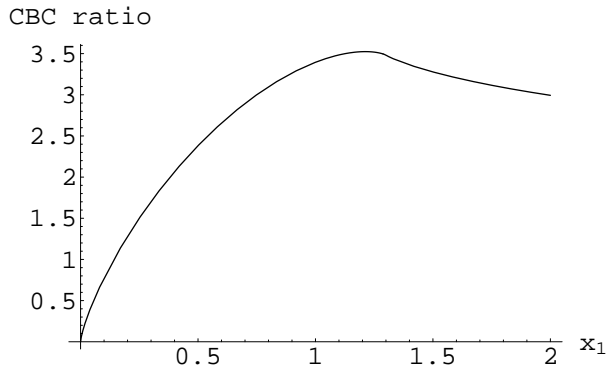}
\caption{\label{fig:Peak1}CBC ratio as a function of the first  turning
point ($x_{1}^{(2)}$). The value of $x^{(1)}_{1}$, which gave 
(Fig.~\ref{fig:quantization}) a ratio of 3.48049,
was 1.30083.
This can be optimally improved to 1,
by choosing $x^{(2)}_{1}=0.141784$.}
\end{figure}

\begin{figure}
\includegraphics{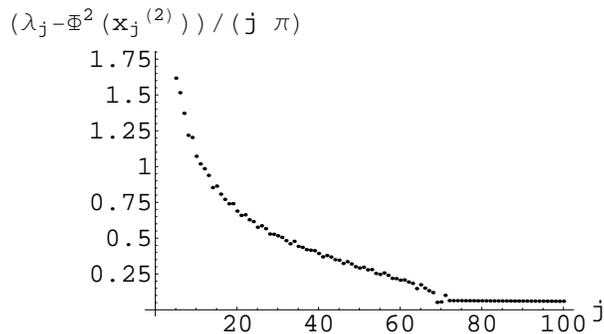}
\caption{\label{fig:RelDev}The difference, divided by $j \pi$, of the
$j$-th Riemann zero minus that value for the zero predicted 
from (\ref{SUSY}) using
the {\it revised} set of turning points ($x_{j}^{(2)}$) 
(Fig.~\ref{fig:Improvement}) 
and the {\it original} 
phases ($\alpha_{k}^{(1)}$) (Fig.~\ref{fig:100})}
\end{figure}

In {\it all} one hundred cases, when we sought to drive the CBC ratios closer to unity, the new 
turning points 
chosen, $x_{j}^{(2)}$, were {\it smaller} than the originals $x_{j}^{(1)}$
(Fig.~\ref{fig:Improvement}). 
(There is some anomalous behavior apparent, in that $x_{71}^{(2)} < 
x_{70}^{(2)}$.)
For the first seventy-one cases, the new turning points ($x_{j}^{(2)}$) do
yield the desired 
CBC ratio of 1 (Fig.~\ref{fig:Peak1}). (We have not been able to obtain a revised set of $\alpha^{(2)}_{k}$'s based on the 100 
 $x^{(2)}_{j}$'s to {\it closely} satisfy (\ref{SUSY}), as we had 
with the 100 $x^{(1)}_{j}$'s.)
When the {\it revised} set of turning points 
$x_{j}^{(2)}$ was employed, along with 
the {\it original} set of phases $\alpha_{k}^{(1)}$'s 
(Fig.~\ref{fig:100}), instead of highly
accurate predictions of the $\lambda_{j}$'s, 
we obtained signficant discrepancies. (The sum-of-squares measure of
fit to (\ref{SUSY}) jumped enormously from $2.68927 \cdot 10^{-14}$ to 129157.)
In Fig.~\ref{fig:RelDev}, we show the differences (all positive) of 
the $\lambda_{j}$'s  minus their predicted values, divided by $j \pi$.

In this paper, we have {\it not} adhered to  the 
apparent misstatement in
\cite[eq. (89)]{carlos1}, that the Riemann zeros ($\lambda_{j}$) should be
``properly normalized'' (``unfolded'' \cite[eq. (6)]{katz}) as
\begin{equation}
\lambda_{j} \rightarrow \frac{\lambda_{j}}{2 \pi} \log{\lambda_{j}},
\end{equation}
but rather conducted our analyses in terms of the ``raw'' 
unnormalized Riemann zeros
themselves ($\lambda_{j}$'s), 
similarly to the work of Wu and Sprung \cite{wu} (cf. \cite{wu2}).
\subsection{Non-iterative (simultaneous) approach}
\subsubsection{$n=m=8$ analyses} \label{nm8}

We commenced an additional series of analyses --- in lieu of the 
iterative scheme just proposed and discussed (but not fully pursued, due to
apparent computational demands).
For the case $n=m=8$, we tried to {\it simultaneously}
estimate the scaled phases ($\alpha_{k}$) and 
turning points ($x_{j}$), {\it as well as} 
the frequency parameter $\gamma$ (for 
a total of seventeen parameters). 
We sought to minimize --- as we will throughout the remainder of 
this study --- a certain sum-of-squares, one of the two addends 
being composed of 
the squared 
differences between the LHS and RHS of the 
(SUSY-QM) equations (\ref{SUSY}) and the 
other addend of the squared differences 
between the LHS and RHS of the (CBC) equations (\ref{CBC}).
(We simply, in effect, 
weighted the two sums {\it equally}, but there certainly appears --- in 
retrospect --- to 
be opportunity for refinement in the choice of weights, since 
the second (CBC) sum seemed to consistently contribute {\it considerably} 
more weight.)
We also --- as mentioned --- allowed the parameter $\gamma$ 
that determines the frequencies $\gamma^k$ to vary (between 1 and 5), 
rather than holding it fixed at 3 (as we had preliminarily thought 
there was some numerical evidence  
for doing so earlier 
(sec.~\ref{Preliminary})).

The minimum sum-of-squares we achieved was 181.6414. A predominant 
part of this 
sum --- 130.236 --- was attributable to the lack of fit to the 
CBC quantization conditions (\ref{CBC}).
{\it Remarkably}, all eight  phases ($\alpha_{k}$) were estimated to be 
essentially {\it zero} (or, effectively, the same $2 \pi$). The {\it largest} 
deviation of these eight
estimates from 0 or $2 \pi$ 
was, in fact, only $8.347 \cdot 10^{-8}$. (We had abandoned also
here our earlier computational 
constraint that the phases be monotonically nondecreasing and only required, 
numerically, that they lie between 0 and $2 \pi$.)
For $\gamma$ we obtained an estimate of 1.4119.
The eight
estimated turning points ($x_{j}$) were 
(0.662734, 1.17022, 1.58516, 2.28689, 2.32938, 3.18343, 3.21295, 3.30561), 
which can be compared with (1.30083, 1.87866, 2.20626, 
2.64243, 2.84142, 3.20489, 3.4613, 3.64459), 
based on the smooth
Wu-Sprung potential (\ref{start}).
Based on these two sets of turning points --- and 
$\alpha_{k}=0$, $k=1,\ldots,8$ and $\gamma=1.41119$ --- the 
corresponding CBC ratios (in the same order of presentation) are 
(1.80006, 1.55891, 1.3248, 1.32142, 1.22064, 1.23564, 1.23128, 1.21166)
and
(1.46964, 1.34903, 1.19796, 1.26407, 1.14572, 1.23298, 1.201, 1.16905). 
We can see that the latter set of CBC ratios is clearly 
superior/more desirable  --- giving 
rise to a sum-of-squares deviation of only 83.5937 {\it vs.} 
the 130.236 mentioned-before --- but 
use of the original set of turning points ($x_{j}^{(1)}$)
yields a much larger 
sum-of-squares for deviations from the other (SUSY-QM) set 
of equations (\ref{SUSY}) of 620.878, as opposed to the $51.4054 
= 181.6414 -130.236$ 
we obtained.

\subsubsection{$n=m=7$ analyses} \label{nm7}
We had also conducted some parallel $n=m=7$ analyses.
We obtained from Mathematica,
three minima (with figures of merit of 69.6285, 140.887 and 
171.939). 
The corresponding estimates of $\gamma$ 
were 1.16585, 1.34436 and 1.4099. 
Of the 21 ($7 \times 3$) phases estimated, none was
greater than 0.000210 removed from either 0 or $2 \pi$.

We will find, apparently importantly, in Sec.~\ref{redo} that for this 
$n=m=7$ case, we can obtain a much {\it improved} figure of merit
of 13.703 by incorporating one additional parameter ($\sigma$) into
the SUSY-QM model (\ref{SUSY}), one which {\it scales} the
Weierstrass fractal function $W(x,\gamma,\frac{3}{2},\alpha_{k})$.
Before doing so, however, we will examine, in detail, the fractal
Wu-Sprung potential.

\section{Analyses of the Wu-Sprung {\it fractal} potential} \label{potential}
\begin{figure}[hb]
\includegraphics{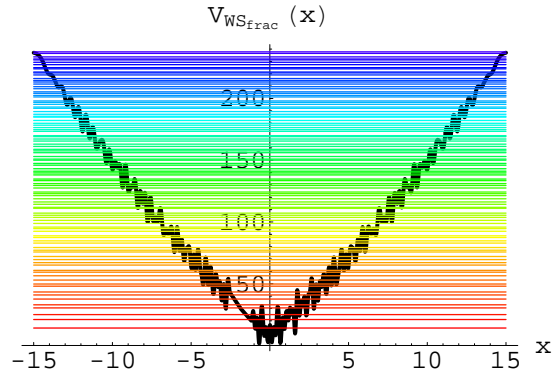}
\caption{\label{fig:True}{\it Fractal} Wu-Sprung potential, $V_{WS_{frac}}(x)$,
 based on the first
100 Riemann zeros}
\end{figure}

\begin{figure}
\includegraphics{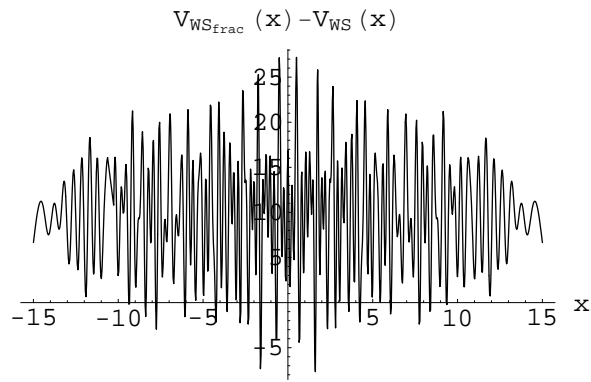}
\caption{\label{fig:Nonsmooth}{\it Non-smooth} component of the Wu-Sprung 
fractal potential $V_{WS_{frac}}(x)$, that is, $V_{WS_{frac}}(x) 
-V_{WS}(x)$.}
\end{figure}

\begin{figure}]
\includegraphics{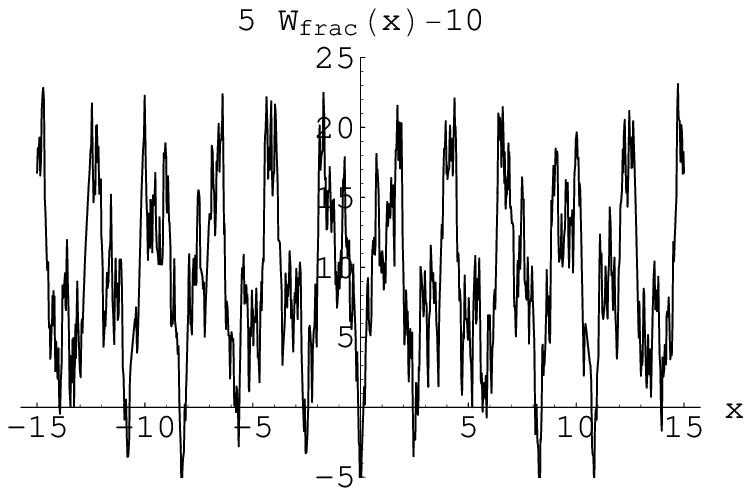}
\caption{\label{fig:WeierAdjust}Weierstrass fractal function contribution
to (\ref{SUSY}) multiplied by 5, then diminished by 10, with $\gamma=2.3$ 
and $\alpha_{k}=0$, with intention of adjusting the Weierstrass fractal 
function to 
resemble Fig.~\ref{fig:Nonsmooth} as closely as possible}
\end{figure}

We implemented, for the case  $n=100$, with an algorithm provided by 
M. Trott, 
the {\it dressing
transformation} for the {\it inversion} of 
the eigenvalues ($\lambda_{j}$'s), 
employed by van Zyl and Hutchinson \cite{vanZyl}, in order 
to obtain the {\it fractal} form $V_{WS_{frac}}(x)$ 
of the Wu-Sprung potential 
(Fig.~\ref{fig:True}) (cf. \cite[Fig. 2]{wu}, \cite[Figs. 1, 3]{vanZyl}).
In Fig.~\ref{fig:Nonsmooth} we show the {\it non-smooth} 
component (the ``residual'' in statistical parlance) of the fractal
WS potential. In other words, we have subtracted away from Fig~\ref{fig:True}, the 
smooth WS potential $V_{WS}(x)$.
We have tried to manipulate the Weierstrass fractal function (\ref{Weierfractal}), as 
incorporated into (\ref{SUSY}) --- assuming all $\alpha_{k}=0$ --- so 
as to closely resemble Fig.~\ref{fig:Nonsmooth}.
A simply visual, non-rigorous analysis lead us to present 
Fig.~\ref{fig:WeierAdjust}, in which we take $\gamma=2.3$ {\it and}
multiply the Weierstrass fractal function term in (\ref{SUSY}) 
by 5 and subtract 10. (Possibly, the CM model might be modified, somewhat
in line with this observation. So, it would appear, that even though the
Weierstrass fractal function, with its fractal dimension $D$ set to 
$\frac{3}{2}$ and the fractal component of the $V_{WS_{frac}}$ share the same
fractal dimension, they seem to differ considerably in certain of their
overall [scale, location] properties. We will explore this issue
in Sec.~\ref{Scale}.)
\subsection{Turning points ($x_{j_{h}}^{frac}$) 
of the {\it fractal} Wu-Sprung potential}
Using the results of the dressing transformation, 
we initially considered taking for the values of 
$x_{j}^{frac}$,
those closest to the original smooth Wu-Sprung turning point 
($x_{j}^{(1)}$), 
for which $\lambda_{j} = V_{WS_{frac}}(x_{j}^{frac})$, as our new/revised 
turning point to employ in the CM model.

Pertaining to the issue, for the edification of the reader, we will note that 
in an e-mail exchange of ours with M. Trott, he 
had commented:  ``But they 
(van Zyl and Hutchinson \cite{vanZyl})
use an exact solution of the Schr\"odinger equation.
So in their approach they never need turning points.
This bring us back to my remark from the beginning of this project.
What is a 'turning point' in a fractal potential? If you take the smallest
$x_{Min}$, such that $\int_{-x_{Min}}^{x_{Min}} \sqrt{e -V(x)}   = n \hbar $
with $e = V(x_{Min})$, then, because of the fractal nature, the barrier
could be very small and the tunnelling contributions from larger
$x$ could be quite significant.'' 
In response, Castro first wrote:`` 
Since the Weierstrass function oscillates, there are many 
possible choices for the $x_{j}^{frac}$, so
you have to choose for the turning points those which are the closest to the ones [$x_{j}^{(1)}$] 
obtained from the smooth Wu and Sprung potential alone. This process will determine the $x_{j}^{frac}$ points uniquely''.
Trott further 
suggested: ``From an 
ordinary quantum mechanics point of view using the 'nearest' ones seems
unconvincing to me. Already for the simple case of two wells separated
by a wall, the classical JWKB formulas get nontrivial corrections'' (cf. \cite{robnik}).
(However, in retrospect, 
it now appears that the most effective manner in which to resolve
this issue --- at least, in pragmatic terms --- is to choose that
$x^{frac}_{j}$ for which the associated fractal CBC relation (\ref{CBC}) is most
closely satisfied.)

\begin{figure}
\includegraphics{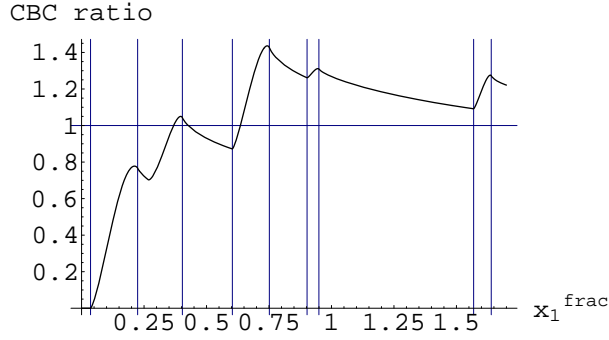}
\caption{\label{fig:CBCFRAC1}CBC ratios based on the fractal Wu-Sprung 
potential $V_{WS_{frac}}(x)$. The nine ($h=1,\ldots,9$) 
turning points ($x_{1_{h}}^{frac}$)
for which $V_{WS_{frac}}(x_{1_{h}}^{frac}) =\lambda_{1}$ are indicated.}
\end{figure}

\begin{figure}
\includegraphics{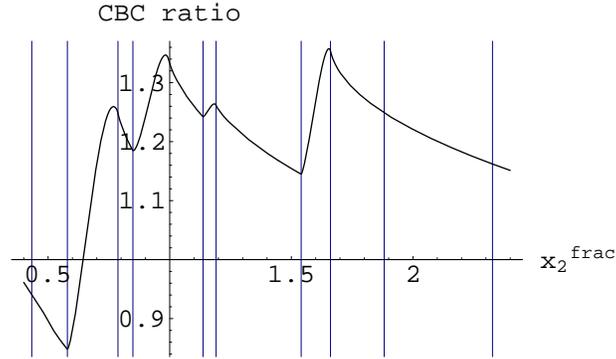}
\caption{\label{fig:CBCFRAC2}CBC ratios based on the fractal Wu-Sprung
potential $V_{WS_{frac}}(x)$. The ten ($h=1,\ldots,10$)
turning points ($x_{2_{h}}^{frac}$)
for which $V_{WS_{frac}}(x_{2_{h}}^{frac}) =\lambda_{2}$ are indicated.}
\end{figure}

\begin{figure}
\includegraphics{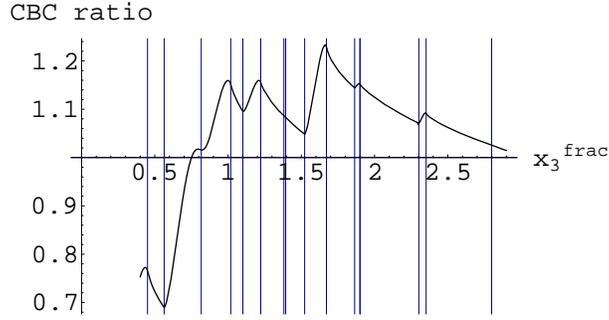}
\caption{\label{fig:CBCFRAC3}CBC ratios based on the fractal Wu-Sprung
potential $V_{WS_{frac}}(x)$. The sixteen ($h=1,\ldots,16$)
turning points ($x_{3_{h}}^{frac}$)
for which $V_{WS_{frac}}(x_{3_{h}}^{frac}) =\lambda_{3}$ are indicated.}
\end{figure}

\begin{figure}
\includegraphics{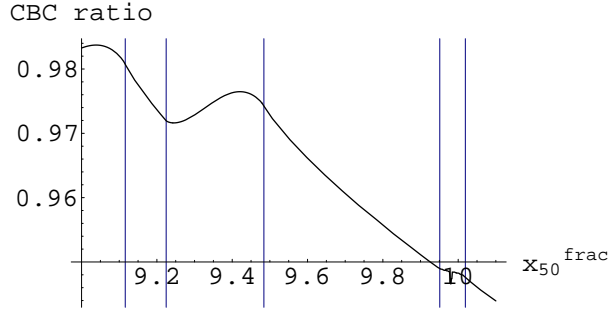}
\caption{\label{fig:CBCFRAC50}CBC ratios based on the fractal Wu-Sprung
potential $V_{WS_{frac}}(x)$. The five ($h=1,\ldots,5$)
turning points ($x_{50_{h}}^{frac}$)
for which $V_{WS_{frac}}(x_{50_{h}}^{frac}) =\lambda_{3}$ are indicated.}
\end{figure}

We, in fact,
found {\it nine} possible values of $x_{1_{h}}^{frac}$, that is
(0.0359371, 0.224349, 0.403363, 
0.60336, 0.75138, 0.902281, 0.949646, 1.5688, 1.63925), satisfying the turning point relation
$V_{WS_{frac}}(x_{1_{h}}^{frac}) =\lambda_{1}$. The correponding CBC
ratios  (Fig.~\ref{fig:CBCFRAC1}) --- (0., 0.766538, 1.04323,
0.871811, 1.42656, 1.26132, 1.30892, 1.09161, 1.27263) seemed,
except obviously for
the first,  reasonably
well-behaved.
For the second Riemann zero ($\lambda_{2}$), the ten turning points 
found were 
(0.432856, 0.57922, 0.787248, 0.848629, 1.13689, 1.19066, 1.54042, 1.66097, \
1.88178, 2.32703) with corresponding 
CBC ratios (Fig.~\ref{fig:CBCFRAC2}) 
(0.966061, 0.848087, 1.2481, 1.18487, 1.24222, 1.26174, 
1.14487, 1.35361, 
1.24891, 1.16204).
 For $\lambda_{3}$, 
the sixteen turning points
(0.450109, 0.563705, 0.816881, 1.0183, 1.10108, 1.22306, 1.38131, 
1.39341, 1.5228, 1.67159, 1.86394, 
1.90039, 1.90039, 2.30223, 2.35184, 2.79924)
with corresponding CBC ratios (Fig.~\ref{fig:CBCFRAC3}) of 
(0.765244, 0.689671, 1.01581, 1.15142, 1.09523, 1.15534, 
1.08533, 1.08337, 
1.04846, 1.22877, 1.14428, 1.15193, 1.15193, 1.07091, 1.0915, 1.02537).
We also show results (Fig.~\ref{fig:CBCFRAC50}) for $x_{50_{h}}^{frac}$,
with turning points (9.11597, 9.22462, 9.4841, 9.9511, 10.019) 
and fractal CBC ratios 
(0.980756, 0.971956, 0.97427, 0.948974, 0.94763).
\subsection{Formulas for certain fractal turning points}
Continuing along these lines,
we found that particular
value (0.949646) of $x$ for which $V_{WS_{frac}}(x) = \lambda_{1}$, and 
which, of the nine solutions,  was, at the same time, {\it smaller} than 
$x_{1}^{(1)} = 1.30083$, and {\it closest} to it. With the
use of ``Plouffe's inverter'' (http://pi.lacim.uqam.ca/eng/), 
we obtained a remarkable
formula that would appear to apply to it:
\begin{equation} \label{x1frac}
x_{1}^{frac} = 0.949646 \approx 0.99999996  \cdot 10^{-6} 
\frac{t e^{\lambda_{1}}}{GS},
\end{equation}
where $t \approx 1.8392867$ is the {\it tribonacci constant} 
and $GS$ is the {\it Gelfond-Schneider} constant, $2^{\sqrt{2}}$.

Proceeding further along these lines, 
the value ($x=1.660974$) we obtained
for $x_{2}^{frac}$, that satisfied $V_{WS_{frac}}(x)= \lambda_{2}$, and was 
smallest and closest to $x_{2}^{(1)}= 1.87866$ was expressible as
\begin{equation} \label{x2frac}
x_{2}^{frac} = 1.660974 \approx 1.0000028005 \cdot 10^{-21} 
\frac{e^{2 \lambda_{3}}}{t \log^2{(2 +\sqrt{3})}}.
\end{equation}

Additionally, 
\begin{equation}
 x_{3}^{frac}= 1.9003895  \approx 0.99999998 \cdot 10^{-40} 
e^{3 \lambda_{5}} 
\frac{3 Trott}{(2 + \sqrt{3})^2},
\end{equation}
where $Trott = 
0.010841015122311136$ is the {\it Trott constant}, and $x_{3}^{(1)}= 2.20626$. 
(The Trott constant has the property that it is invariant, in a certain 
sense. If one expands the number
digit by digit as if it were a continued fraction, then
the number remains the same \cite{finch}.)

Further, 
\begin{equation}
x_{4}^{frac}= 2.3843247 \approx 1.00000003037 \cdot 
10^{-72} \frac{e^{4 \lambda_{7}}}{t^2 Cahen^2 \log{\zeta{(5)}}},
\end{equation}
where $Cahen \approx 0.6434105462883$ is the {\it Cahen constant}, used in the 
theory of continued fractions. ($x_{4}^{(1)} = 2.64243$ and $\zeta$ 
is the Riemann zeta function --- {\it exact} formulas for which are known 
for {\it even} integral arguments.)

Also,
\begin{equation}
x_{5}^{frac} = 2.8338417 \approx 1.00000004769 \cdot 10^{-106} 
\frac{\zeta{(3)}  e^{4 {\sqrt{2}}} e^{5 \lambda_{9}}}{\Gamma{(\frac{7}{12})}}.
\end{equation}
($x_{5}^{frac}$ is rather close to $x_{5}^{(1)} = 2.84142$. $\zeta(3)$ is 
known to be {\it irrational} and is sometimes referred to as 
``Ap\'ery's constant'' \cite{apery}.)
\section{{\it Scaling} the Weierstrass Fractal Function} \label{Scale}
The analyses underlying Figs.~\ref{fig:Nonsmooth} and \ref{fig:WeierAdjust}, 
and their comparison,  
indicated to us that
possibly the CM model of the Riemann zeros might be enhanced by
multiplying the Weierstrass fractal contribution (\ref{SUSY}) by some 
(new) scaling 
factor ($\sigma$). (Since in Sec.~\ref{Preliminary}, we were essentially 
finding parameters to set the Weierstrass function to zero, the issue
of scaling was not germane there.) So, we modified the type of analyses 
reported in Secs.~\ref{nm8} and \ref{nm7} to estimate $\sigma$, as well as the 
(scaled) phases ($\alpha_{k}$), the turning points ($x_{j}$) and 
the frequency parameter ($\gamma$).
\subsection{$n=m=7$ analyses} \label{redo}
We use exactly 
the same criterion of best (unscaled) 
fit, employed previously, that is the total
sum-of-squares between the LHS's and RHS's of (\ref{SUSY}) and (\ref{CBC}).
The sum-of-squares we obtained was 13.703 (8.48114 coming from
the lack of fit to the CBC 
quantization conditions (\ref{CBC})), {\it much} superior 
to the best fit of 69.6285 in Sec.~\ref{nm7}. 
This, of course, lends strong support for the relevance of incorporating
the scaling parameter $\sigma$, which was set to
3.92036, along with $\gamma = 2.18081$ --- {\it both} 
estimates being 
quite consistent with
Fig.~\ref{fig:WeierAdjust}. 
The seven estimated phases were
(0.915274, 6.28319, 6.28319, 1.20429, 5.33637, 0.700917, 0.0) --- three of
them being essentially 0 or $2 \pi \approx 6.28318531$.
The seven estimated turning points were 
(0.321253, 0.676572, 0.936234, 1.65921, 1.79613, 2.1688, 2.18378), 
considerably smaller than in Secs.~\ref{nm7} and \ref{nm8}.
The associated CBC ratios --- all now relatively close 
to unity --- are (1.1386, 1.16732, 0.951104, 1.13392, 
1.02881, 1.067, 1.06951).

At this stage, we could have continued this $n=m=7$ analysis to see if
we could obtain further minima with figures of merit less than 13.703.
But now knowing that Mathematica could produce results for $n=m=7$, we 
were interested in seeing how far we could ``push'' Mathematica for cases
$n=m>7$. So, rather than now dwelling on a particular $n=m$ case,
we kept seeking higher-order  instances. (At this point, we treated each
analysis independently, that is we did not employ results of
earlier [lower-dimensional] analyses as initial guesses.)
We did obtain single minima for $n=m=8$ and $n=m=9$, but the figures of
merit were rather weak (39.1008 and 85.689), so we do not detail the
results here.

\subsection{$n=m=10$ analyses} \label{redo10}
Now we expand to the case of ten phases and ten turning points.
The (relatively quite {\it excellent}) 
sum-of-squares we obtained was 2.86609 (2.68872 coming from
the lack of fit to the CBC 
quantization conditions (\ref{CBC})).
The estimate of $\gamma$ was 1.30466 and of $\sigma$, 
1.49575.
The ten estimated phases were 
(6.28319, 6.28319, 1.27142, 0.0493542, 6.28319, 0.0555025, 
0.0000178313, 6.28319, 0., $6.00153 \cdot 10^{-23}$),
{\it most} of them being essentially 0 or $2 \pi \approx 6.28318531$.
The ten estimated turning points were 
(0.38402, 0.665384, 0.915631, 1.32033, 1.42245, 
1.83676, 1.95767, 2.33568, 2.68742, 3.06041),
considerably smaller than in Secs.~\ref{nm7} and \ref{nm8}, where
the scaling parameter ($\sigma$) was absent.
The associated CBC ratios --- most now relatively close 
to unity --- are (1.27306, 1.142, 1.00359, 1.03838, 
0.969627, 1.02007, 1.0058, 0.978438, 1.00323, 0.984601).

So, an important unresolved question is whether or not
ideally {\it all} the phases ($\alpha_{k}$) should be set to zero, or
whether the fact that all the phases were estimated as zero in 
Sec.~\ref{nm8}
was principally a manifestation of the need for a scaling
parameter {\it greater} than unity. 
This is because if one sets all the phases to zero, the
term-by-term interference in the Weierstrass fractal function
is minimized and its overall contribution/influence/magnitude increased. 
(I thank Carlos Castro for this last observation.)
\subsubsection{Zero phases} \label{Zero}
To address this last-discussed issue, 
still within the $n=m=10$ scenario, we {\it a priori} set all ten phases
($\alpha_{k}$) to zero, leaving us with twelve parameters (ten turning points, 
plus $\gamma$ and $\sigma$) to estimate.
We obtained four minima, one relatively inferior (24.0743) and three, all
yielding approximately 7.2358 (6.758 of this from the CBC equations) --- but all more than  the 2.86609 obtained immediately above --- for 
the sum-of-squares, and all with
essentially the same set of estimates.
The estimates were $\gamma=1.25274$, $\sigma=1.23289$ (both interestingly close to one another) and for the turning
points (0.406011, 0.72083, 0.88735, 1.21503, 1.38215, 1.72066, 
2.07865, 2.24067, 2.81678, 2.91862).
The CBC ratios were
(1.4716, 1.20207, 1.04624, 1.0322, 0.953197, 0.985181, 0.987027, 0.9602, 
1.03323, 0.995425).

So, no yet fully convincing evidence for the 
proposition that {\it all} the phases ($\alpha_{k}$) 
should be set to zero, has appeared in this $n=m=10$ suite of analyses
(though Castro believes this is an interesting hypothesis well worth pursuing). 
\subsection{$n=m=12$}
We were able to obtain two minima ({\it independently}-generated
with different random seeds) for this case, having associated
sums-of-squares of 4.7436 (4.51953 attributable to the fit 
to the CBC quantization conditions) and 4.88101 (4.59739 
so-attributable).

These two minima had rather similiar estimates of
$\gamma$, $\sigma$, as well as of the twelve turning points.
(For the smaller of the two, $\gamma$ and $\sigma$ were estimated
as 1.16902 and 1.90025, and for the larer, 1.17209 and 1.79825.)
The turning points for the smaller (of the two small) minima were:
(0.370638, 0.518695, 0.65557, 1.3734, 1.49501, 
1.90065, 2.13024, 2.25476, 2.83058, 2.96561, 3.18434, 3.69962),
and for the larger of the two 
(0.371207, 0.526829, 0.697696, 1.33586, 1.45628, 
1.91474, 2.13455, 2.25897, 2.76173, 2.88164, 3.10636, 3.50429).
The phases ($2 \pi \alpha_{k}$) for the smaller of the two 
were: (5.06388, 1.19176, 5.10283, 5.10239, 1.54197, 4.87505, 
5.75582, 5.07329, 6.28319, 6.28319, 0.491806, 0.) (three being essentially
0 or $2 \pi$) and 
(1.36199, 6.28319, 5.32029, 3.14014, 6.28319, 0.932408, 
1.48183, 6.2407, 0.811536, 6.28186, 0.0160344, 0.)
(five now being essentially or quite close to 0 or $2 \pi$).
The CBC ratios based on the smallest (of the three $n=m=10$) 
minima were: (1.45746, 1.09132, 0.891351, 1.03399, 0.967626, 1.0088, 
1.01311, 0.981312, 1.01104, 0.986514, 1.00474, 1.00054).
\subsection{$n=m=15$} \label{redo15}
Mathematica yielded a minimum having a sum-of-squares of 10.8781,
10.1905 of the total
stemming from deviations from the idealized CBC relationships.

The estimate of $\gamma$ was 1.10999 and of 
the scaling parameter $\sigma$, 2.39239.
The fifteen phases were estimated to be 
(5.85115, 1.90972, 6.07707, 3.38707, 5.79917, 
3.01828, 5.21112, 0.972974, 6.28319, 4.06528, 5.13804, 
1.09178, 5.42993, 0.0000400221, 0.887475) and the fifteen turning points, 
(0.427238, 0.575873, 0.674962, 0.964599, 1.67025, 
1.86543, 1.99707, 2.09904, 
2.87102, 2.96694, 3.12156, 3.26476, 3.37853, 3.42232, 3.66246).
The resultant CBC ratios were (1.62295, 1.205, 0.970261, 0.919054, 
0.964421, 1.03902, 1.01992, 0.974149, 
1.00741, 0.984101, 1.00373, 1.01898, 1.01317, 0.977452, 0.996592).

\begin{figure}
\includegraphics{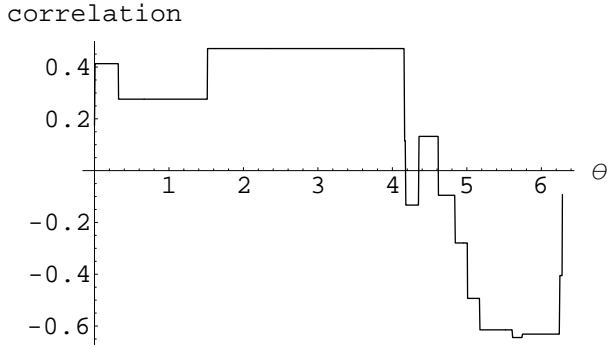}
\caption{\label{fig:Correlation}Correlation coefficients between the first
set of phases for the $n=m=15$ case, and the second set of phases, with
the addition to its members of an angle $\theta \in [0, 2 \pi]$. 
The highest plateau 
of level 0.471313  extends 
from $\theta= 1.56108995$ to 4.1583035}
\end{figure}

Subsequently, we obtained a second minimum (starting with different 
random seeds), with a slightly inferior
sum-of-squares of 11.6345 (10.9979 due to failure to completely
satisfy the CBC relations). The estimate (1.10638) 
of $\gamma$ obtained was notably 
similar to the first estimate (1.10999), while the estimate of $\sigma$ 
was 1.92895 {\it vs.} 2.39239, previously. 
The estimated phases were
(1.44074, 6.28318, 1.66813, 0.542437, 2.10471, 4.76708, 2.12488, 
0.0381873, 1.27475, 0.672527, 1.10953, 1.92716, 5.96083, 0.00013761, 0.) 
and the estimated turning points,
(0.447782, 0.605638, 0.711853, 1.03604, 1.65431, 1.87966, 2.02289, 
2.13206, 2.85323, 2.95601, 3.11569, 3.25825, 3.37216, 3.41606, 3.70782).
The {\it correlation coefficient} between the two sets of turning points
was extremely high, 0.999796, while that between the two sets of 
phases was positive, but apparently quite weak, 0.174294.
(The sets of turning points are constrained to be strictly increasing,
while the phases are not.)
This can be (maximally) improved to 0.471313 by adding $\theta = \pi$ 
to the second set of phases
(Fig.~\ref{fig:Correlation}).
The correlation between the two sets of CBC ratios was extremely high, 
0.997186. (Obviously, no ordering constraints has been placed upon them.)

\subsection{$n=m=20$} \label{redo20}
Initially, we obtained a minimum  with a somewhat disappointingly large
(in comparison with previous analyses) sum-of-squares of 
60.1355 (57.0834 stemming from the failure to fully satisfy the CBC 
relations). The estimate of $\gamma$ was 1.03784.
Then, we obtained a somewhat superior sum-of-squares of 
53.6492 (50.1484 coming from the lack of complete fit to the CBC relations).
The estimate of $\gamma$, 1.049,  was quite similar.

Ultimately though, we found a much superior solution with a sum-of-squares
equal to 11.4018, of which 10.3952 was attributable to the failure to fully
fit the CBC relations. The estimate of $\gamma$ was 1.08285 and of
$\sigma$, 1.23124.
The twenty phases were 
(6.28319, 5.12918, 1.82172, 0.627236, 2.28354, 6.28057,
  2.09377, 2.27574, 6.28319, 6.28319, 1.39899, 1.05003, 4.08259,
  0.00871371, 6.28319, 5.20743,
 0.00527643, 6.28319, $3.08389 \cdot 10^{-7}$, 0.00224054) --- half of them being quite proximate 
to 0 or $2 \pi$ --- and the twenty turning points, 
(0.430261, 0.581042, 0.682297, 0.990577, 1.66733, 1.8869,
  2.03975, 2.16318, 2.7836, 2.8999, 3.11796, 3.36483, 3.61917,
  3.71294, 4.1214, 4.23216, 4.39019, 4.59062, 5.39774,
  5.49934). 
(Possibly, the relatively large number of phases near either 0 or $2 \pi$ 
is a computational artifact, since these are the {\it extreme} values 
between which the 
phases are constrained to lie.)
Additionally, the resultant CBC ratios were
(1.63002, 1.21105, 0.97568, 0.925948, 0.959536, 1.03755, 1.02111, 0.978159, \
1.01279, 0.983753, 0.998834, 1.01199, 1.00774, 0.97948, 1.01524, 1.00143, \
0.997587, 0.993286, 1.00448, 0.994628).

We have been investigating the $n=m=25$ case (both with unconstrained phases 
and phases set to zero) --- with some programming
improvements provided by M. Trott --- but do not have any specific results
to report presently. (It appears that several days on a PowerMac are
required to locate a single relative  minimum.)
\subsubsection{Decreasing trend in estimates of $\gamma$} \label{gamma}
We see an undeniable trend in the several analyses above. As
the number of unknowns ($n+m+2$) increases, the parameter $\gamma$
{\it decreases}, possibly indicative of a convergence to 
its theoretical lower bound for the Weierstrass fractal function of 1.
(For $n=m=7$, our one estimate of $\gamma$ was 2.18081; for 
$n=m=10$, the estimate was 1.30466; for $n=m=12$,  we had two estimates 
of 1.16902 and 1.17209; for $n=m=15$, 1.10999 and 1.10638; and
for $n=m=20$, 1.03784, 1.049 and 1.08285.)
Szulga and Molz have shown that in the limit $\gamma \to 1$
(with random phases, in the case $D=\frac{3}{2}$ before us here),  ``the Mandelbrot-Weierstrass process is a complex fractional Brownian motion'' \cite{szulga} 
(cf. \cite{szulga2,pipiras,rawlings,biane}). In the Castro-Mahecha model, one takes
the {\it real} part of the Mandelbrot-Weierstrass function (which CM term the
Weierstrass function). To statistically 
test whether or not the estimated series of phases we obtain is
{\it random} in nature, one can apply ``Rao's spacing test'' \cite{russell}. 
(Also, the Rayleigh test and the Kuiper's V test are sometimes employed.)
Applying this test to the first 
two sets of twenty phases reported for our
$n=m=20$ analyses in Sec.~\ref{redo20}, we obtained statistics of 
192.199 and 158.86, respectively.
Any result greater than 192.17 is statistically significant
at the $p=0.001$ level of probability, and any 
result greater than 154.31, at the $p=0.10$ level 
\cite[Table II]{russell}. There is, thus, rather strong evidence that 
the hypothesis that 
the phases are random should be {\it rejected}.

\section{Discussion} \label{discussion}
On the scientifically important question of possible {\it falsifiability}, 
the Castro-Mahecha supersymmetric model of the nontrivial Riemann zeros --- and that of Wu and Sprung \cite{wu} --- would {\it fail} (as Carlos Castro 
indicated) if it could be 
demonstrated 
that the shape of the potential were {\it multifractal}, so that the 
fractal dimension would no longer be a {\it constant} ($D=\frac{3}{2}$).
We are not cognizant of any 
evidence of this, though. However, it might be of some 
value to formally test for such a possibility. 
(Though we are aware of no multifractal counterpart to the Weierstrass
function, Castro suggested that one could sum (or possibly integrate) over 
the fractal dimension ($D$).) Additionally, still within a ``unifractal''
framework, Castro put forth the idea of integrating 
(thus, generalizing) the Weierstrass
fractal function --- with fractal dimension $D$ held constant at 
$\frac{3}{2}$ --- over the parameter 
$\gamma$, using some
appropriate weighting function $f(\gamma)$.

\begin{figure}
\includegraphics{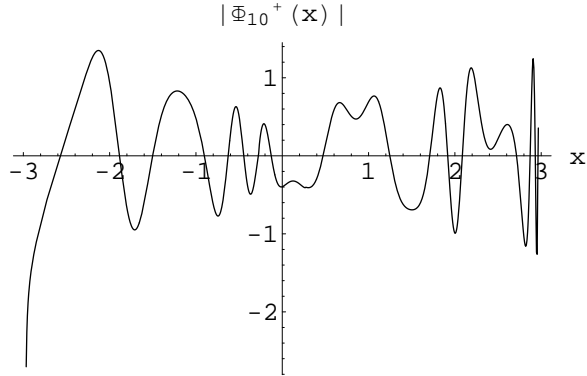}
\caption{\label{fig:Wavefunction}Absolute value ($|\Phi_{10}^{+}(x)|$)
of the wavefunction corresponding to the tenth 
eigenvalue $\lambda_{10}$ and the 
tenth turning point $x_{10}=2.96694$ reported for the $n=m=15$ analysis
(smaller of the two minima) 
in Sec.~\ref{redo15}}
\end{figure}

The presumptive Hilbert-Polya operator would have eigenvalues
equal to the $\lambda_{j}$'s. 
It might be of interest to attempt to estimate its eigen{\it functions}
within the same 
postulated quasi-classical (CBC) framework employed in the CM model,
and in our analyses above. In Sec.~6.3 of \cite{junkerbook} 
(also \cite[eq. (56)]{inomata}) there are presented formulas
for the quasi-classical wave functions. 
The appropriate fractal extensions of these 
formulas --- following the {\it ansatz} employed 
in (\ref{CBC}) --- would need to be employed (cf. \cite{karol}). 
An initial effort of ours along these lines ---  by simple way of 
a single illustration --- yielded 
Fig.~\ref{fig:Wavefunction}.
In the context of their foundational article 
\cite[p. 2597]{wu} 
(cf. \cite{tomiya}), ``Riemann zeros and a fractal potential'', Wu and Sprung noted
that ``All the wave functions decay exponentially beyond the turning point, 
and thus quickly die away. This implies that the lower energy levels have 
very small influence on the potential beyong their turning points \ldots 
Clearly, the potential in the low energy range has greater 
'responsibility' than that in a higher energy range, thus it has more
structure \ldots Based on this argument, the finest structure will be 
determined by the wave number of the last wave function.''

It would be of interest to estimate {\it all} the $n$ wavefunctions in our
various $n=m$ analyses above, and examine whether or not the fractal SUSY
Schr\"odinger equation (\ref{difficult}) is (at least, 
approximately) satisfied.
This would be, in some sense, a test of the appropriateness of the
fractal extension of the CBC relations postulated by Castro and Mahecha
(adopted by them 
in order to avoid having to directly solve (\ref{difficult})).
(It would also, of course, be of interest to 
attempt to solve (\ref{difficult}) more directly, rather than having
to resort to the fractal CBC {\it ansatz} of CM.)

\begin{figure}
\includegraphics{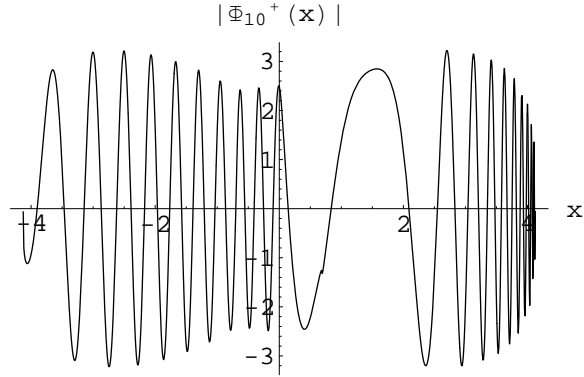}
\caption{\label{fig:WaveFunction15NoFractal}Counterpart to Fig.~\ref{fig:Wavefunction} based simply
on the {\it smooth} Wu-Sprung potential, {\it without} any
fractal contribution, and original (non-fractal) formulas for the
quasi-classical eigenfunctions}
\end{figure}

In Fig.~\ref{fig:WaveFunction15NoFractal}
we show the counterpart to Fig.~\ref{fig:Wavefunction} based simply
on the {\it smooth} Wu-Sprung potential, without any
fractal contribution, and original (non-fractal) formulas 
(\cite[eq. (6.40)]{junkerbook}, \cite[eq. (56)]{inomata}) for the
quasi-classical eigenfunctions.
It would be an interesting exercise to further duplicate some of the analyses
conducted here, but with the replacement of the fractal CBC equations
(\ref{CBC}) with the (conventional) CBC equations.

Throughout the minimization 
analyses above, for the sake of simplicity (and to avoid the 
possible confusion of issues), we have relied upon a quite elementary
least-squares fitting. Modifying this approach to give more equable
weighting to the two sets of coupled equations, and perhaps to lower and 
higher Riemann zeros and turning points might prove useful.

CM believed ``that their fractal SUSY QM model, once the optimum value 
for the amplitude factor $\gamma$ is known, has a great chance of truly
reproducing the zeta zeros, and proving the [Riemann Hypothesis], 
by simply establishing a one-to-one correspondence among the values
of the infinite phases of our Weierstrass function with the zeta zeros'' 
 \cite[sec.~6]{carlos1}. We have tried to make a contribution in such a 
direction here.

Castro suggested that the most
salient/interesting aspect of our several analyses above had emerged in 
Sec.~\ref{Zero}, where for the case $n=m=10$ we had {\it set} the ten
phases to zero (doing so which had been strongly suggested
by the results 
of the immediately preceding
analysis in which the phases were {\it not} constrained, and allowed to vary
between 0 and $2 \pi$, but found, {\it nevertheless}, to be predominantly 
0 or $2 \pi$). 
The estimates of the frequency parameter 
$\gamma$ and the scaling parameter $\sigma$ were quite 
intriguingly close
(1.25274 and 1.23289).
However, the figure of merit, 7.2358, 
in this zero-phase analysis was somewhat
inferior to that [2.86609] in that unconstrained 
$n=m=10$ analysis, which 
rather induced
us here not to subsequently pursue solely further analyses ($n=m>10$) 
in which the phases had
been {\it a priori} set to zero --- but rather let them 
freely vary. Also, we suspected, as previously mentioned, 
that the predominance of 0's or $2 \pi$'s might be, in some way, a 
computational artifact, since these are the imposed limits 
on the phases. (The minimization procedure, we speculated,  
might drive the phases to these
end points, beyond which it could not further proceed 
to seek improvements, and would therefore {\it terminate} there. 
Castro suggested that setting the range of possible phases 
to be $[-\pi,\pi]$ would obviate this 
particular problem --- but possibly 
introduce new ones. We intend to investigate such a direction.) 

Perhaps somewhat relatedly, as 
one of their ``spectral speculations'', Berry and Keating 
\cite[p. 260]{berrykeating}
remarked that the ``Maslov phases associated with the [classical periodic
orbits of the Riemann dynamics] are also peculiar: they are all $\pi$.
The result appears paradoxical in view of the relations between these
phases and the winding numbers of the stable and unstable manifolds
associated with periodic orbits , but finds an explanation in a scheme of
Connes''. (The Maslov phases appear when one
modifies the Bohr-Sommerfield quantization rules by equating the WKB 
integrals (orbits) to $j \pi$ plus a Maslov phase factor, instead of  $j \pi$.)

Castro suggested that the possibility  the Weierstrass fractal function
might assume the form of classical Brownian motion was an unappealing one,
seeing that this would apparently imply that the phases ($\alpha_{k}$) must be random.
(In Sec.~\ref{gamma} we deduced certain 
statistical evidence that the estimated
phases [allowed to freely vary] were {\it not}, in fact, random. 
Of course, we know this to be strictly true, since the phases are constructed
by our estimation procedure.)

\begin{figure}
\includegraphics{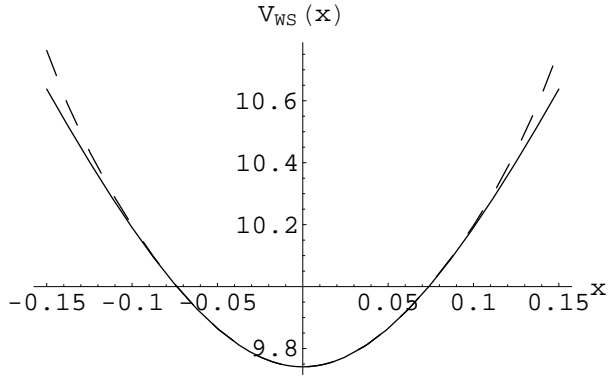}
\caption{\label{fig:DD}Inverted smooth Wu-Sprung 
potential (\ref{WSpotential}) --- using the Mathematica Interpolation
command on a data table of resolution $\frac{1}{40}$ --- and 
the approximation (dashed lines) to the inverted potential 
using the first three 
terms (\ref{first3}) of the Dominici series expansion
(\ref{DDseries})}
\end{figure}

\hspace{2.6in} {\bf{Appendix}}
\newline
At our suggestion, D. Dominici investigated 
the problem of inverting the 
Wu-Sprung potential (\ref{WSpotential}) \cite[eq. (7)]{wu}. 
His analyses  have 
been presented in a ``preliminary report'' \cite{diego}.
Dominici  
gives explicitly 
the first {\it ten} coefficients ($a_{k}$) of a power series expansion,
\begin{equation} \label{DDseries}
V_{WS}(x)=V_{0} + \Sigma_{k=1}^{\infty} a_{k} (\pi x)^{2 k} \omega^{2 k -1} 
(-V_{0})^{1-k},
\end{equation}
where
\begin{equation}
\omega= \Big[\ln{\Big(\frac{V_{0}}{2 \pi} \Big)} \Big]^{-1}.
\end{equation}
(For the critical value $V_{0} = 2 \pi$, above which $V_{WS}(x)$ is 
{\it single}-valued, 
the expansion (\ref{DDseries}) 
is no longer valid. 
A separate analysis \cite[Sec.~2.2]{diego} is, then, of interest.)
Using just the first three of these coefficients,
\begin{equation} \label{first3}
a_{1}= \omega; a_{2} =\frac{4}{3} \omega^2; a_{3} = \frac{8}{15} 
\omega^2 + \frac{28}{9} \omega^3,
\end{equation}
we obtain Fig.~\ref{fig:DD}.

In the analyses reported in the body of this paper, 
we had relied upon the Interpolation
command of Mathematica for the purpose of inversion, but we anticipate
employing the results of Dominici \cite{diego} in our further work.
C. Krattenthaler has indicated that for the 
purpose of studying the {\it asymptotic}  
behavior of the coefficients of the power series, 
it would be best to employ the methodology of Lagrange inversion.
\begin{acknowledgments}
Carlos Castro kindly introduced me to his joint work 
\cite{carlos1} with 
Jorge Mahecha, proposed avenues of 
research to pursue, and discussed in great detail 
the analytical developments.
Further, I would like 
to express gratitude to the Kavli Institute for Theoretical
Physics (KITP)
for computational support,  
to Michael Trott
of Wolfram Research Inc. for his more than generous assistance and expertise
in terms of Mathematica computations, as well as for 
his critical interest and comments, and to Diego Dominici 
for his interest in the problem of inverting the Wu-Sprung potential 
\cite{diego}.
\end{acknowledgments}

\bibliography{Carlos5}

\begin{thebibliography}{32}
\expandafter\ifx\csname natexlab\endcsname\relax\def\natexlab#1{#1}\fi
\expandafter\ifx\csname bibnamefont\endcsname\relax
  \def\bibnamefont#1{#1}\fi
\expandafter\ifx\csname bibfnamefont\endcsname\relax
  \def\bibfnamefont#1{#1}\fi
\expandafter\ifx\csname citenamefont\endcsname\relax
  \def\citenamefont#1{#1}\fi
\expandafter\ifx\csname url\endcsname\relax
  \def\url#1{\texttt{#1}}\fi
\expandafter\ifx\csname urlprefix\endcsname\relax\def\urlprefix{URL }\fi
\providecommand{\bibinfo}[2]{#2}
\providecommand{\eprint}[2][]{\url{#2}}

\bibitem[{\citenamefont{Conrey}(2003)}]{conrey}
\bibinfo{author}{\bibfnamefont{J.~B.} \bibnamefont{Conrey}},
  \bibinfo{journal}{Notices Amer. Math. Soc.} \textbf{\bibinfo{volume}{50}},
  \bibinfo{pages}{341} (\bibinfo{year}{2003}).

\bibitem[{\citenamefont{Bunimovich and Dettmann}(2005)}]{circular}
\bibinfo{author}{\bibfnamefont{L.~A.} \bibnamefont{Bunimovich}}
  \bibnamefont{and} \bibinfo{author}{\bibfnamefont{C.~P.}
  \bibnamefont{Dettmann}}, \bibinfo{journal}{Phys. Rev. Lett.}
  \textbf{\bibinfo{volume}{94}}, \bibinfo{pages}{100201}
  (\bibinfo{year}{2005}).

\bibitem[{\citenamefont{Rockmore}(2005)}]{rockmore}
\bibinfo{author}{\bibfnamefont{D.}~\bibnamefont{Rockmore}},
  \emph{\bibinfo{title}{Stalking the Riemann Hypothesis}}
  (\bibinfo{publisher}{Pantheon}, \bibinfo{address}{New York},
  \bibinfo{year}{2005}).

\bibitem[{\citenamefont{Castro and Mahecha}(2004)}]{carlos1}
\bibinfo{author}{\bibfnamefont{C.}~\bibnamefont{Castro}} \bibnamefont{and}
  \bibinfo{author}{\bibfnamefont{J.}~\bibnamefont{Mahecha}},
  \bibinfo{journal}{Int. J. Geom. Meth. Mod. Phys.}
  \textbf{\bibinfo{volume}{1}}, \bibinfo{pages}{751} (\bibinfo{year}{2004}).

\bibitem[{\citenamefont{Wu and Sprung}(1993)}]{wu}
\bibinfo{author}{\bibfnamefont{H.}~\bibnamefont{Wu}} \bibnamefont{and}
  \bibinfo{author}{\bibfnamefont{W.~L.} \bibnamefont{Sprung}},
  \bibinfo{journal}{Phys. Rev. E} \textbf{\bibinfo{volume}{48}},
  \bibinfo{pages}{2595} (\bibinfo{year}{1993}).

\bibitem[{\citenamefont{Berry and Lewis}(1980)}]{berry}
\bibinfo{author}{\bibfnamefont{M.~V.} \bibnamefont{Berry}} \bibnamefont{and}
  \bibinfo{author}{\bibfnamefont{Z.~V.} \bibnamefont{Lewis}},
  \bibinfo{journal}{Proc. Roy. Soc. London} \textbf{\bibinfo{volume}{370}},
  \bibinfo{pages}{459} (\bibinfo{year}{1980}).

\bibitem[{\citenamefont{Khuri}(2002)}]{khuri}
\bibinfo{author}{\bibfnamefont{N.~N.} \bibnamefont{Khuri}},
  \bibinfo{journal}{Math. Phys. Anal. Geom.} \textbf{\bibinfo{volume}{5}},
  \bibinfo{pages}{1} (\bibinfo{year}{2002}).

\bibitem[{\citenamefont{Rosu}(2003)}]{rosu}
\bibinfo{author}{\bibfnamefont{H.~C.} \bibnamefont{Rosu}},
  \bibinfo{journal}{Mod. Phys. Lett. A} \textbf{\bibinfo{volume}{18}},
  \bibinfo{pages}{1205} (\bibinfo{year}{2003}).

\bibitem[{\citenamefont{van Zyl and Hutchinson}(2003)}]{vanZyl}
\bibinfo{author}{\bibfnamefont{B.~P.} \bibnamefont{van Zyl}} \bibnamefont{and}
  \bibinfo{author}{\bibfnamefont{D.~A.~W.} \bibnamefont{Hutchinson}},
  \bibinfo{journal}{Phys. Rev. E} \textbf{\bibinfo{volume}{67}},
  \bibinfo{pages}{066211} (\bibinfo{year}{2003}).

\bibitem[{\citenamefont{Laskin}(2002)}]{laskin}
\bibinfo{author}{\bibfnamefont{N.}~\bibnamefont{Laskin}},
  \bibinfo{journal}{Phys. Rev. E} \textbf{\bibinfo{volume}{66}},
  \bibinfo{pages}{056108} (\bibinfo{year}{2002}).

\bibitem[{\citenamefont{Inomata and Junker}(1994)}]{inomata}
\bibinfo{author}{\bibfnamefont{A.}~\bibnamefont{Inomata}} \bibnamefont{and}
  \bibinfo{author}{\bibfnamefont{G.}~\bibnamefont{Junker}},
  \bibinfo{journal}{Phys. Rev. A} \textbf{\bibinfo{volume}{50}},
  \bibinfo{pages}{3638} (\bibinfo{year}{1994}).

\bibitem[{\citenamefont{Edwards}(1974)}]{edwards}
\bibinfo{author}{\bibfnamefont{H.~M.} \bibnamefont{Edwards}},
  \emph{\bibinfo{title}{Riemann's zeta function}} (\bibinfo{publisher}{Academic
  Press}, \bibinfo{address}{New York}, \bibinfo{year}{1974}).

\bibitem[{\citenamefont{Ma and Hori}(2003)}]{mahori}
\bibinfo{author}{\bibfnamefont{C.}~\bibnamefont{Ma}} \bibnamefont{and}
  \bibinfo{author}{\bibfnamefont{Y.}~\bibnamefont{Hori}}, in
  \emph{\bibinfo{booktitle}{ASME 2003 Design Engineering Technical
  Conferences}} (\bibinfo{year}{2003}), \bibinfo{note}{(DET2003VIB-48736)}.

\bibitem[{\citenamefont{Dominici}()}]{diego}
\bibinfo{author}{\bibfnamefont{D.}~\bibnamefont{Dominici}},
  \eprint{math.CA/0510341}.

\bibitem[{\citenamefont{Bhaduri et~al.}(2005)\citenamefont{Bhaduri, Sakhr,
  Sprung, Dutt, and Suzuki}}]{bhaduri}
\bibinfo{author}{\bibfnamefont{R.~K.} \bibnamefont{Bhaduri}},
  \bibinfo{author}{\bibfnamefont{J.}~\bibnamefont{Sakhr}},
  \bibinfo{author}{\bibfnamefont{D.~W.~L.} \bibnamefont{Sprung}},
  \bibinfo{author}{\bibfnamefont{R.}~\bibnamefont{Dutt}}, \bibnamefont{and}
  \bibinfo{author}{\bibfnamefont{A.}~\bibnamefont{Suzuki}},
  \bibinfo{journal}{J. Phys. A} \textbf{\bibinfo{volume}{38}},
  \bibinfo{pages}{L183} (\bibinfo{year}{2005}).

\bibitem[{\citenamefont{Planat and Rosu}(2003)}]{planat}
\bibinfo{author}{\bibfnamefont{M.}~\bibnamefont{Planat}} \bibnamefont{and}
  \bibinfo{author}{\bibfnamefont{H.~C.} \bibnamefont{Rosu}},
  \bibinfo{journal}{Phys. Lett. A} \textbf{\bibinfo{volume}{315}},
  \bibinfo{pages}{1} (\bibinfo{year}{2003}).

\bibitem[{\citenamefont{R{\"u}schendorf}(1995)}]{ludger}
\bibinfo{author}{\bibfnamefont{L.}~\bibnamefont{R{\"u}schendorf}},
  \bibinfo{journal}{Ann. Statist.} \textbf{\bibinfo{volume}{23}},
  \bibinfo{pages}{1160} (\bibinfo{year}{1995}).

\bibitem[{\citenamefont{Katz and Sarnak}(1999)}]{katz}
\bibinfo{author}{\bibfnamefont{N.~M.} \bibnamefont{Katz}} \bibnamefont{and}
  \bibinfo{author}{\bibfnamefont{P.}~\bibnamefont{Sarnak}},
  \bibinfo{journal}{Bull. Amer. Math. Soc.} \textbf{\bibinfo{volume}{36}},
  \bibinfo{pages}{1} (\bibinfo{year}{1999}).

\bibitem[{\citenamefont{Wu et~al.}(1990)\citenamefont{Wu, Valli{\'e}res, Feng,
  and Sprung}}]{wu2}
\bibinfo{author}{\bibfnamefont{H.}~\bibnamefont{Wu}},
  \bibinfo{author}{\bibfnamefont{M.}~\bibnamefont{Valli{\'e}res}},
  \bibinfo{author}{\bibfnamefont{D.~H.} \bibnamefont{Feng}}, \bibnamefont{and}
  \bibinfo{author}{\bibfnamefont{D.~W.~L.} \bibnamefont{Sprung}},
  \bibinfo{journal}{Phys. Rev. A} \textbf{\bibinfo{volume}{42}},
  \bibinfo{pages}{1027} (\bibinfo{year}{1990}).

\bibitem[{\citenamefont{Robnik et~al.}(1999)\citenamefont{Robnik, Salasnich,
  and Vrani{\v c}ar}}]{robnik}
\bibinfo{author}{\bibfnamefont{M.}~\bibnamefont{Robnik}},
  \bibinfo{author}{\bibfnamefont{L.}~\bibnamefont{Salasnich}},
  \bibnamefont{and} \bibinfo{author}{\bibfnamefont{M.}~\bibnamefont{Vrani{\v
  c}ar}}, \bibinfo{journal}{Nonlin. Phen. Complex Syst.}
  \textbf{\bibinfo{volume}{2}}, \bibinfo{pages}{49} (\bibinfo{year}{1999}).

\bibitem[{\citenamefont{Finch}(2003)}]{finch}
\bibinfo{author}{\bibfnamefont{S.~R.} \bibnamefont{Finch}},
  \emph{\bibinfo{title}{Mathematical Constants}} (\bibinfo{publisher}{Cambridge
  U. P.}, \bibinfo{address}{Cambridge}, \bibinfo{year}{2003}).

\bibitem[{\citenamefont{Ap{\'e}ry}(1979)}]{apery}
\bibinfo{author}{\bibfnamefont{R.}~\bibnamefont{Ap{\'e}ry}},
  \bibinfo{journal}{Ast{\'e}rique} \textbf{\bibinfo{volume}{61}},
  \bibinfo{pages}{11} (\bibinfo{year}{1979}).

\bibitem[{\citenamefont{Szulga and Molz}(2001)}]{szulga}
\bibinfo{author}{\bibfnamefont{J.}~\bibnamefont{Szulga}} \bibnamefont{and}
  \bibinfo{author}{\bibfnamefont{F.}~\bibnamefont{Molz}}, \bibinfo{journal}{J.
  Statist. Phys.} \textbf{\bibinfo{volume}{104}}, \bibinfo{pages}{1317}
  (\bibinfo{year}{2001}).

\bibitem[{\citenamefont{Szulga}(2002)}]{szulga2}
\bibinfo{author}{\bibfnamefont{J.}~\bibnamefont{Szulga}},
  \bibinfo{journal}{Statist. Prob. Lett.} \textbf{\bibinfo{volume}{56}},
  \bibinfo{pages}{301} (\bibinfo{year}{2002}).

\bibitem[{\citenamefont{Pipiras and Taqqu}(2000)}]{pipiras}
\bibinfo{author}{\bibfnamefont{V.}~\bibnamefont{Pipiras}} \bibnamefont{and}
  \bibinfo{author}{\bibfnamefont{M.~S.} \bibnamefont{Taqqu}},
  \bibinfo{journal}{Fractals} \textbf{\bibinfo{volume}{8}},
  \bibinfo{pages}{369} (\bibinfo{year}{2000}).

\bibitem[{\citenamefont{Rawlings}(2003)}]{rawlings}
\bibinfo{author}{\bibfnamefont{P.~K.} \bibnamefont{Rawlings}},
  \bibinfo{journal}{J. Statist. Phys.} \textbf{\bibinfo{volume}{111}},
  \bibinfo{pages}{769} (\bibinfo{year}{2003}).

\bibitem[{\citenamefont{Biane et~al.}(2001)\citenamefont{Biane, Pitman, and
  Yor}}]{biane}
\bibinfo{author}{\bibfnamefont{P.}~\bibnamefont{Biane}},
  \bibinfo{author}{\bibfnamefont{J.}~\bibnamefont{Pitman}}, \bibnamefont{and}
  \bibinfo{author}{\bibfnamefont{M.}~\bibnamefont{Yor}},
  \bibinfo{journal}{Bull. Amer. Math. Soc.} \textbf{\bibinfo{volume}{38}},
  \bibinfo{pages}{435} (\bibinfo{year}{2001}).

\bibitem[{\citenamefont{Russell and Levitin}(1995)}]{russell}
\bibinfo{author}{\bibfnamefont{G.~S.} \bibnamefont{Russell}} \bibnamefont{and}
  \bibinfo{author}{\bibfnamefont{D.~J.} \bibnamefont{Levitin}},
  \bibinfo{journal}{Commun. Statist.: Simul. Comput.}
  \textbf{\bibinfo{volume}{24}}, \bibinfo{pages}{879} (\bibinfo{year}{1995}).

\bibitem[{\citenamefont{Junker}(1996)}]{junkerbook}
\bibinfo{author}{\bibfnamefont{G.}~\bibnamefont{Junker}},
  \emph{\bibinfo{title}{Supersymmetric Methods in Quantum and Statistical
  Physics}} (\bibinfo{publisher}{Springer-Verlag}, \bibinfo{address}{Berlin},
  \bibinfo{year}{1996}).

\bibitem[{\citenamefont{W{\'o}jcik et~al.}(2000)\citenamefont{W{\'o}jcik,
  Bia{\l}ynicki, and {\.Z}yczkowski}}]{karol}
\bibinfo{author}{\bibfnamefont{D.}~\bibnamefont{W{\'o}jcik}},
  \bibinfo{author}{\bibfnamefont{I.}~\bibnamefont{Bia{\l}ynicki}},
  \bibnamefont{and}
  \bibinfo{author}{\bibfnamefont{K.}~\bibnamefont{{\.Z}yczkowski}},
  \bibinfo{journal}{Phys. Rev. Lett.} \textbf{\bibinfo{volume}{85}},
  \bibinfo{pages}{5022} (\bibinfo{year}{2000}).

\bibitem[{\citenamefont{Tomiya and Sakamoto}(2003)}]{tomiya}
\bibinfo{author}{\bibfnamefont{M.}~\bibnamefont{Tomiya}} \bibnamefont{and}
  \bibinfo{author}{\bibfnamefont{S.}~\bibnamefont{Sakamoto}},
  \bibinfo{journal}{e-J. Surf. Sci. Nanotech.} \textbf{\bibinfo{volume}{1}},
  \bibinfo{pages}{175} (\bibinfo{year}{2003}).

\bibitem[{\citenamefont{Berry and Keating}(1999)}]{berrykeating}
\bibinfo{author}{\bibfnamefont{M.~V.} \bibnamefont{Berry}} \bibnamefont{and}
  \bibinfo{author}{\bibfnamefont{J.~P.} \bibnamefont{Keating}},
  \bibinfo{journal}{SIAM Review} \textbf{\bibinfo{volume}{41}},
  \bibinfo{pages}{236} (\bibinfo{year}{1999}).

\end{thebibliography}

\end{document}